\documentclass[10pt]{article}
\newtheorem{thm}{Theorem}

\newtheorem{cor}[thm]{Corollary}

\newcommand{\Z}{\hbox{\bf Z}}
\newcommand{\G}{\mbox{$\Gamma$}}

\newcommand{\beeq}{\begin{eqnarray*}}
\newcommand{\eneq}{\end{eqnarray*}}
\newcommand{\proof}{\noindent {\it Proof.\hspace{4mm}}}

\newcommand{\qfd}{\hfill $\fbox{}$\vspace{4mm}}
\def\newpic#1{%
\def\emline##1##2##3##4##5##6{%
\put(##1,##2){\special{em:point #1##3}}%
\put(##4,##5){\special{em:point #1##6}}%
\special{em:line #1##3,#1##6}}}
\newpic{}
\def\emline#1#2#3#4#5#6{%
\put(#1,#2){\special{em:moveto}}%
\put(#4,#5){\special{em:lineto}}}
\def\newpic#1{}
\title{From the Coxeter graph to the Klein graph}
\author{Italo J. Dejter
\\ University of Puerto Rico \\ Rio Piedras, PR 00936-8377 \\ ijdejter@uprrp.edu}
\date{}

\begin{document}
\maketitle

\begin{abstract} We show that the 56-vertex Klein cubic graph $\G'$ can be obtained from the 28-vertex Coxeter cubic graph $\G$ by 'zipping' adequately the squares of the 24 7-cycles of $\G$ endowed with an orientation obtained by considering $\G$ as a $\mathcal C$-ul\-tra\-ho\-mo\-ge\-neous digraph, where $\mathcal C$ is the collection formed by both the oriented 7-cycles $\vec{C}_7$ and the $2$-arcs $\vec{P}_3$ that tightly fasten those $\vec{C}_7$ in $\G$. In the process, it is seen that $\G'$ is a ${\mathcal C}'$-ultrahomogeneous (undirected) graph, where ${\mathcal C}'$ is the collection formed by both the 7-cycles $C_7$ and the $1$-paths $P_2$ that tightly fasten those $C_7$ in $\G'$. This yields an embedding of $\G'$ into a 3-torus $T_3$ which forms the Klein map of Coxeter notation $(7,3)_8$. The dual graph of $\G'$ in $T_3$ is the distance-regular Klein quartic graph, with corresponding dual map of Coxeter notation $(3,7)_8$.
\end{abstract}

\noindent{\bf Keywords:} ultrahomogeneous graph; digraph; shortest cycle; arc-transitivity

\noindent{\bf 2000 Mathematics subject classification:} 05C62, 05B30, 05C20, 05C38

\section{Introduction}

\noindent The study of ultrahomogeneous graphs (resp. digraphs) can be traced back to \cite{Sheh,Gard,Ronse,Cam,GK},
(resp. \cite{Fra,Lach,Cher}). Following a line of research initiated in \cite{I},  given a collection $\mathcal C$ of (di)graphs closed under isomorphisms, a (di)graph $G$ is said to be $\mathcal C$-{\it ul\-tra\-ho\-mo\-ge\-neous} (or $\mathcal C$-UH) if every isomorphism between two induced members of $\mathcal C$ in $G$ extends to an auto\-mor\-phism of $G$. If ${\mathcal C}=\{H\}$ is the isomorphism class of a (di)graph $H$, we say that such a $G$ is $\{H\}$-UH or $H$-UH.
In \cite{I}, $\mathcal C$-UH graphs are defined and studied when ${\mathcal C}$ is the collection of either {\bf(a)} the complete graphs, or {\bf(b)} the disjoint unions of complete graphs, or {\bf(c)} the complements of those unions.\bigskip

\noindent We may consider a graph $G$ as a digraph by considering  each edge $e$ of $G$ as a pair of oppositely oriented (or O-O) arcs $\vec{e}$ and $(\vec{e})^{-1}$. Then, {\it zipping} or {\it fastening} $\vec{e}$ and $(\vec{e})^{-1}$, operation that we define as uniting $\vec{e}$ and $(\vec{e})^{-1}$, allows to obtain precisely $e$, a simple technique to be used below.
(In \cite{Dx}, however, a strongly connected $C_4$-UH oriented graph without O-O arcs was presented). In other words,
$G$ must be a graph considered as a digraph, that is, for any two vertices $u,v\in V(G)$, the arcs $\vec{e}=(u,v)$ and $(\vec{e})^{-1}=(v,u)$ are both present in the set $A(G)$ of arcs of $G$, with the union $\vec{e}\,\cup(\vec{e})^{-1}$ interpreted as the (undirected) edge $e\in E(G)$ of $G$. If we write $\vec{f}=(\vec{e})^{-1}$, then clearly $(\vec{f})^{-1}=\vec{e}$ and $f=e$.\bigskip

\noindent Let $M$ be a subgraph of a graph $H$ and let $G$ be both an $M$-UH and an $H$-UH graph.
We say that $G$ is an $\{H\}_{M}$-{\it UH  graph} if, for each copy $H_0$ of $H$ in $G$ containing a copy $M_0$ of $M$, there exists exactly one copy $H_1\neq H_0$ of $H$ in $G$ with $V(H_0)\cap V(H_1)=V(M_0)$ and $E(H_0)\cap E(H_1)=E(M_0)$.
The vertex and edge conditions above can be condensed as $H_0\cap H_1=M_0$. We say that such a $G$ is {\it tightly fastened}. This is generalized by saying that an $\{H\}_M$-UH graph $G$ is an $\ell${\it -fastened} $\{H\}_M$-UH graph if given a copy $H_0$ of $H$ in $G$ containing a copy $M_0$ of $M$, then there exist exactly $\ell$ copies $H_i\neq H_0$ of $H$ in $G$ such that $H_i\cap H_0\supseteq M_0$, for each $i=1,2,\ldots,\ell$, with at least $H_1\cap H_0=M_0$.\bigskip

\noindent Now, let $\vec{M}$ be a subdigraph of a digraph $\vec{H}$ and let the graph $G$ be both an $\vec{M}$-UH and an $\vec{H}$-UH digraph. We say that $G$ is an $\{\vec{H}\}_{\vec{M}}$-{\it UH digraph} if for each copy $\vec{H}_0$ of $\vec{H}$ in $\vec{G}$ containing a copy $\vec{M}_0$ of $\vec{M}$ there exists exactly one copy $\vec{H}_1\neq\vec{H}_0$ of $\vec{H}$ in $G$ with $V(\vec{H}_0)\cap V(\vec{H}_1)=V(\vec{M}_0)$ and $A(\vec{H}_0)\cap\bar{A}(\vec{H}_1)=A(\vec{M}_0)$, where $\bar{A}(\vec{H}_1)$ is formed by those arcs $(\vec{e})^{-1}$ whose orientations are reversed with respect to the orientations of the arcs $\vec{e}$ of $A(\vec{H}_1)$. Again, we may say that such a $G$ is {\it tightly fastened}. This case is used in the construction of Section 3.\bigskip

\noindent The Coxeter graph $\G=F_{028}A$ \cite{F} is a distance-transitive hypohamiltonian \cite{Bondy} cubic graph of order $n=28$, diameter $d=4$, girth $g=7$, arc-transitivity $k=3$, having exactly $\eta=24$ $g$-cycles, $a=336$ automorphisms, intersection array ${\mathcal I}=\{3,2,2,1;1,1,1,2\}$ and weakly regular parameters ${\mathcal W}=(28,(3),(0),(0,1))$.
The Klein cubic graph $\G'=F_{056}B$ is a hamiltonian cubic graph with $n'=2n$, $d'=6$, $g'=g$, $k'=2$, $\eta'=\eta$, $a'=a$ and ${\mathcal W}'=(24,(7),(2),(0,2))$, (not to be confused with the bipartite double graph of $\G$, denoted $F_{056}C$); see \cite{F,Sch,K,SL}. (We remark that $\G$ can be obtained as the graph whose vertices are the 6-cycles of the Heawood graph $\G''=F_{014}A$ \cite{F}, with any two vertices adjacent if and only if the 6-cycles they represent are disjoint, where we recall that $\G''$ is a distance-transitive hamiltonian cubic graph with $n''=14$, $d''=3$, $g''=6$, $k''=4$, $\eta''=n$, $a''=a$, ${\mathcal I}''=\{3,2,2;1,1,3\}$ and ${\mathcal W}''=(14,(3),(0),(0,1))$.)\bigskip

\noindent Given a finite graph $H$ and a subgraph $M$ of $H$ with $|V(H)|>3$, we say that a graph $G$ is ({\it strongly fastened}) or {\it SF} $\{(H\}_M$-{\it UH} if there is a sequence of connected subgraphs $M=M_1,M_2\ldots,M_t\equiv K_2$ such that: {\bf(a)} $M_{i+1}$ is obtained from $M_i$ by the deletion of a vertex, for $i=1,\ldots,t-1$ and {\bf(b)} $G$ is a $(2^i-1)$-fastened $\{H\}_{M_i}$-UH graph, for $i=1,\ldots,t$.
Theorem 1 below asserts that $\G$ is an SF $\{C_7\}_{P_3}$-UH graph.\bigskip

\noindent Theorem 2 establishes that $\G$ is a $\{\vec{C_7}\}_{\vec{P_3}}$-UH digraph. In Section 3, squaring the resulting oriented 7-cycles allows the recovery of $\G'$ dressed up as a $\{C_7\}_{P_2}$-UH graph, via zipping of the O-O induced 2-arcs shared (as 2-paths) by the pairs of O-O 7-cycles.\bigskip

\noindent As in \cite{F,Sch,K,SL}, the dual graph of $\G'$ with respect to an embedding of its 24 7-cycles into a 3-torus (known as the Klein map, of Coxeter notation $(7,3)_8$, see argument previous to Theorem 3, below) is the Klein quartic graph ${\mathcal K}$ (of Corollary 4), a 24-vertex distance-regular graph with intersection array $\{7,4,1;1,2,7\}$ and weakly regular parameters $(24,(7),(2),(0,2))$.

\section{$\{C_7\}_{P_3}$-UH and $\{\vec{C_7}\}_{\vec{P_3}}$-UH properties of $\G$}

\begin{thm}
$\G$ is an SF $\{C_g\}_{P_3}$-UH graph. In par\-ti\-cu\-lar, $\G$
has exactly $6ng^{-1}=24$ $g$-cycles.
\end{thm}

\proof We have to see that $\G$ is a $(2^{i+1}-1)$-fastened $\{C_g\}_{P_{3-i}}$-UH graph, for $i=0,1$. In fact,
each $(2-i)$-path $P=P_{3-i}$ of $\G$ is shared exactly by $2^{i+1}$ $g$-cycles of $\G$, for $i=0,1$.
This and a simple counting argument for the number of $g$-cycles yield the assertions in the statement.\qfd

\noindent In fact, the proof above can be extended in order to establish that every distance-transitive cubic graph $G$ with girth $=g$ and AT $=k$, (including $G=\G''$), is an SF $\{C_g\}_{P_{i+2}}$-UH graph, for $i=0,1,\ldots,k-2$, and in particular a $\{C_g\}_{P_k}$-UH graph with exactly $2^{k-2}3ng^{-1}$ $g$-cycles.\bigskip

\noindent Given a $\{\vec{C}_g\}_{\vec{P}_k}$-UH graph $G$, an assignment of an orientation to each $g$-cycle of $G$ such that the two $g$-cycles shared by each $(k-1)$-path receive opposite orientations yields a $\{\vec{C}_g\}_{\vec{P}_k}$-{\it O-O assignment}, (or $\{\vec{C}_g\}_{\vec{P}_k}$-OOA). The collection of $\eta$ oriented $g$-cycles corresponding to the $\eta$ $g$-cycles of $G$, for a particular $\{\vec{C}_g\}_{\vec{P}_k}$-OOA will be called an $\{\eta\vec{C}_g\}_{\vec{P}_k}$-OOC. Each such cycle will be expressed with its successive composing vertices expressed between parentheses but without separating commas, (as is the case for arcs $(u,v)$ and 2-arcs $(u,v,w)$), where as usual the vertex that succeeds the last vertex of the cycle is its first vertex.\bigskip

\begin{thm} $\G$ is $\{\vec{C_g}\}_{\vec{P_k}}$-UH, or $\{\vec{C_7}\}_{\vec{P_3}}$-UH.
\end{thm}

\proof $\G$ is obtained from three 7-cycles $(u_1u_2u_3u_4u_5u_6u_0)$, $(v_4v_6v_1v_3v_5v_0v_2)$, $(t_3t_6$ $t_2t_5t_1t_4t_0)$ by adding a copy of $K_{1,3}$ with degree-1 vertices $u_x,v_x,t_x$ and a central degree-3 vertex $z_x$, for each $x\in\Z_7$. Then $\G$ admits the $\{24\,\vec{C}_7\}_{\vec{P}_3}$-OOC:
$$\begin{array}{lll}\vspace*{.5mm}
^{\{\underline{0^1}=(u_1u_2u_3u_4u_5u_6u_0),}_{\,\,\,\underline{1^1}=(u_1z_1v_1\,v_3\,z_3u_3u_2),} &
^{\underline{0^2}=(v_1v_3v_5v_0v_2v_4v_6),}_{\underline{1^2}=(z_4v_4v_2v_0z_0\,t_0t_4),} &
^{\underline{0^3}=(t_1t_5t_2\,t_6\,t_3\,t_0\,t_4),}_{\underline{1^3}=(t_6t_2t_5z_5u_5u_6z_6),}\\
\vspace*{.5mm} ^{\,\,\,\underline{2^1}=(v_5z_5u_5u_4u_3\,z_3\,v_3),}_{\,\,\,\underline{3^1}=(v_5v_0z_0u_0u_6\,u_5\,z_5),} &
^{\underline{2^2}=(t_6z_6v_6v_4v_2\,z_2t_2),}_{\underline{3^2}=(z_4t_4t_1z_1v_1\,v_6v_4),} &
^{\underline{2^3}=(u_1z_1t_1t_4t_0z_0u_0),}_{\underline{3^3}=(t_6t_2z_2u_2u_3z_3t_3),} \\
\vspace*{.5mm}^{\,\,\,\underline{4^1}=(u_1u_0z_0v_0v_2\,z_2\,u_2),}_{\,\,\,\underline{5^1}=(z_4u_4u_3u_2z_2\,v_2\,v_4),} &
^{\underline{4^2}=(t_6t_3z_3v_3v_1v_6\,z_6),}_{\underline{5^2}=(v_5v_3v_1z_1t_1t_5\,z_5),} &
^{\underline{4^3}=(z_4u_4u_5z_5t_5t_1t_4),}_{\underline{5^3}=(t_6z_6u_6u_0z_0t_0t_3),} \\
^{\,\,\,\underline{6^1}=(z_4v_4v_6z_6u_6\,u_5\,u_4),}_{\,\,\,\underline{7^1}=(u_1u_0u_6z_6v_6\,v_1\,z_1),} &
^{\underline{6^2}=(v_5v_3z_3t_3t_0\,z_0v_0),}_{\underline{7^2}=(v_5z_5t_5t_2z_2\,v_2v_0),} &
^{\underline{6^3}=(u_1u_2z_2t_2t_5t_1z_1),}_{\underline{7^3}=(z_4t_4t_0t_3z_3u_3u_4)\}.}
\end{array}$$
\qfd

\noindent In fact, Theorem 2 can be adapted to a statement for every distance-transitive cubic graph which is neither $\G''$ nor the Petersen, Pappus or Foster graphs.

\section{`Zipping' the squares $(\vec{C}_7)^2$ in $\G$ towards $\G'$}

\noindent In this section, we keep using the construction and notation of $\G$ and of its $\{24\vec{C}_7\}_{\vec{P}_3}$-OOC, as conceived in the proof of Theorem 2.
Consider the collection $(\vec{\,\mathcal C}_7)^2(\G)$ of squares \cite{FH} of oriented $7$-cycles in the
$\{24\vec{C}_7\}_{\vec{P}_3}$-OOC of $\G$ in that proof.
From now on, each initial vertex $w_0$ of an arc $\vec{e}=(w_0,w_1)$ of a member $\vec{C}_7^2$ of $(\vec{\,\mathcal C}_7)^2(\G)$, the arc $\vec{e}$ itself and its terminal vertex $w_1$ are respectively indicated by, or marked with, the symbols $v_0,v_1,v_2$ representing the respective vertices of the 2-arc $\vec{E}=(v_0,v_1,v_2)$ of $\vec{C}_7$ associated with $\vec{e}$. For example, if $\vec{C}_7=\underline{0^1}=(u_1u_2u_3u_4u_5u_6u_0)$, so that $\vec{C}_7^2=(\underline{0^1})^2=(u_1u_3u_5u_0u_2u_4u_6)$, then the arc $(u_1,u_3)$ of $\vec{C}_7^2=(\underline{0^1})^2$ is indicated by means of $u_2$, while $u_1$ and $u_3$ are indicated exactly by means of those same symbols, namely $u_1$ and $u_3$. In sum, we get the following indications over $\vec{C}_7^2=(\underline{0^1})^2$:
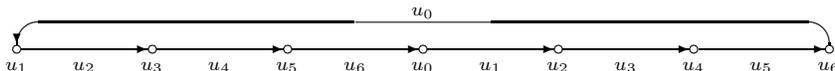
\begin{figure}[htp]
\unitlength=0.60mm
\special{em:linewidth 0.4pt}
\linethickness{0.4pt}
\begin{picture}(191.00,13.00)
\put(10.00,5.00){\circle{2.00}}
\put(40.00,5.00){\circle{2.00}}
\put(10.00,1.00){\makebox(0,0)[cc]{$_{u_1}$}}
\put(25.00,1.00){\makebox(0,0)[cc]{$_{u_2}$}}
\put(40.00,1.00){\makebox(0,0)[cc]{$_{u_3}$}}
\put(70.00,5.00){\circle{2.00}}
\put(55.00,1.00){\makebox(0,0)[cc]{$_{u_4}$}}
\put(70.00,1.00){\makebox(0,0)[cc]{$_{u_5}$}}
\put(100.00,5.00){\circle{2.00}}
\put(85.00,1.00){\makebox(0,0)[cc]{$_{u_6}$}}
\put(100.00,1.00){\makebox(0,0)[cc]{$_{u_0}$}}
\put(130.00,5.00){\circle{2.00}}
\put(115.00,1.00){\makebox(0,0)[cc]{$_{u_1}$}}
\put(130.00,1.00){\makebox(0,0)[cc]{$_{u_2}$}}
\put(160.00,5.00){\circle{2.00}}
\put(145.00,1.00){\makebox(0,0)[cc]{$_{u_3}$}}
\put(160.00,1.00){\makebox(0,0)[cc]{$_{u_4}$}}
\put(190.00,5.00){\circle{2.00}}
\put(175.00,1.00){\makebox(0,0)[cc]{$_{u_5}$}}
\put(190.00,1.00){\makebox(0,0)[cc]{$_{u_6}$}}
\put(84.50,6.00){\oval(149.00,10.00)[lt]}
\emline{85.00}{11.00}{1}{115.00}{11.00}{2}
\put(115.00,6.00){\oval(150.00,10.00)[rt]}
\put(100.00,13.00){\makebox(0,0)[cc]{$_{u_0}$}}
\put(11.00,5.00){\vector(1,0){28.00}}
\put(41.00,5.00){\vector(1,0){28.00}}
\put(71.00,5.00){\vector(1,0){28.00}}
\put(101.00,5.00){\vector(1,0){28.00}}
\put(131.00,5.00){\vector(1,0){28.00}}
\put(161.00,5.00){\vector(1,0){28.00}}
\put(10.00,7.00){\vector(0,-1){1.00}}
\end{picture}
\end{figure}

\noindent where the leftmost horizontal edge stands for the exemplified arc $(u_1,u_3)$. We zip now corresponding O-O arc pairs in the squares $\vec{C}_7^2$ obtained from $\G$, in order to recover $\G'$ with the desired $\mathcal C$-UH properties.
The following sequence of operations is performed:
$$\G\,\,\,\rightarrow\,\,\,\{24\vec{C}_7\}_{\vec{P}_3}\mbox{-OOC}(\G)\,\,\,\rightarrow\,\,\,(\vec{\,\mathcal C}_7)^2(\G)^2(\G)\,\,\,\rightarrow\,\,\, \G'.$$ Next, we explain how this operation $\G\rightarrow \G'$ is composed.
The Fano plane $\mathcal F$, with point set $J_7=\{1,\ldots,7\}$ and line
set $\{124$, $235$, $346$, $457$, $561$, $672$, $713\}$, yields a coloring of the vertices and edges of $\G$, as represented on the upper left quarter of Figure 1, below,
where the color of each vertex $v$ of $\G$ (written in boldface in the next paragraph, for clarity) and the colors of its three incident edges form a quadruple $q$ whose complement $\mathcal F\setminus q$ is used to denote $v$, (\cite{GR} page 69).
Moreover: {\bf(a)} the triple formed by the colors of the edges incident to each $v$ of $\G$ is a line of $\mathcal F$; {\bf(b)} the color of each edge $e$ of $\G$ together with the colors of the endvertices of $e$ form a line of $\mathcal F$.\bigskip

\noindent The vertices $u_x,z_x,v_x,t_x$ created in the presentation of the $\{24\,\vec{C}_7\}_{\vec{P}_3}$-OOC in the proof of Theorem 2 are depicted concentrically in the mentioned representation of $\G$ in Figure 1, from the outside in, starting say  downward from top with colors $x=\mathbf1,\mathbf5,\mathbf4,\mathbf3$ for respective vertices $257={\mathcal F}\setminus\mathbf1364$, $134={\mathcal F}\setminus\mathbf5602$, $567={\mathcal F}\setminus\mathbf4013$, $356={\mathcal F}\setminus\mathbf3214$, which are shown solid in the figure against a backdrop of the remaining hollow vertices.\bigskip

\noindent The squares $\vec{C}_7^2$ corresponding to the 24 oriented 7-cycles $\vec{C}_7$ of $\G$ are represented: {\bf(a)} via their induced cyclically-presented orientations and {\bf(b)} with each vertex $v$ (resp. arc $\vec{e}$) of a $\vec{C}_7^2$ conveniently indicated by means of a color $c(v)$ for $v$ (resp., conveniently indicated by means of a subindex, color $_{c(u)}$ for the middle vertex $u$ of the 2-path $\vec{E}$ of $\vec{C}_7$ that $\vec{e}$ represents). The net effect that this color notation produces makes the 24 oriented 7-cycles $\vec{C}_7^2$ pairwise distinguishable, thus providing them with a distinctive and well-defined presentation. As an example, we go back to the oriented 7-cycle $\vec{C}_7^2=(\underline{0^1})^2$ pictured above, showing now how it receives its colors $c(u_i)$:

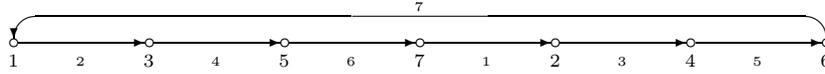
\begin{figure}[htp]
\unitlength=0.60mm
\special{em:linewidth 0.4pt}
\linethickness{0.4pt}
\begin{picture}(191.00,13.00)
\put(10.00,5.00){\circle{2.00}}
\put(40.00,5.00){\circle{2.00}}
\put(10.00,1.00){\makebox(0,0)[cc]{$_1$}}
\put(25.00,1.00){\makebox(0,0)[cc]{$_{_2}$}}
\put(40.00,1.00){\makebox(0,0)[cc]{$_3$}}
\put(70.00,5.00){\circle{2.00}}
\put(55.00,1.00){\makebox(0,0)[cc]{$_{_4}$}}
\put(70.00,1.00){\makebox(0,0)[cc]{$_5$}}
\put(100.00,5.00){\circle{2.00}}
\put(85.00,1.00){\makebox(0,0)[cc]{$_{_6}$}}
\put(100.00,1.00){\makebox(0,0)[cc]{$_7$}}
\put(130.00,5.00){\circle{2.00}}
\put(115.00,1.00){\makebox(0,0)[cc]{$_{_1}$}}
\put(130.00,1.00){\makebox(0,0)[cc]{$_2$}}
\put(160.00,5.00){\circle{2.00}}
\put(145.00,1.00){\makebox(0,0)[cc]{$_{_3}$}}
\put(160.00,1.00){\makebox(0,0)[cc]{$_4$}}
\put(190.00,5.00){\circle{2.00}}
\put(175.00,1.00){\makebox(0,0)[cc]{$_{_5}$}}
\put(190.00,1.00){\makebox(0,0)[cc]{$_6$}}
\put(84.50,6.00){\oval(149.00,10.00)[lt]}
\emline{85.00}{11.00}{1}{115.00}{11.00}{2}
\put(115.00,6.00){\oval(150.00,10.00)[rt]}
\put(100.00,13.00){\makebox(0,0)[cc]{$_{_7}$}}
\put(11.00,5.00){\vector(1,0){28.00}}
\put(41.00,5.00){\vector(1,0){28.00}}
\put(71.00,5.00){\vector(1,0){28.00}}
\put(101.00,5.00){\vector(1,0){28.00}}
\put(131.00,5.00){\vector(1,0){28.00}}
\put(161.00,5.00){\vector(1,0){28.00}}
\put(10.00,7.00){\vector(0,-1){1.00}}
\end{picture}
\end{figure}

\noindent which can be written in short as $(1_23_45_67_12_34_56_7)$, meaning that $c(u_0)=7$, $c(u_i)=i$, for $i=1,\ldots,6$ and if $\vec{e}=(u_i,u_{i+2})$, with $i+j$ taken mod 7 for $j=1,2$, where 0 is rewritten as 7, then $c(\vec{e})=i+1$, this color written as a subindex.\bigskip

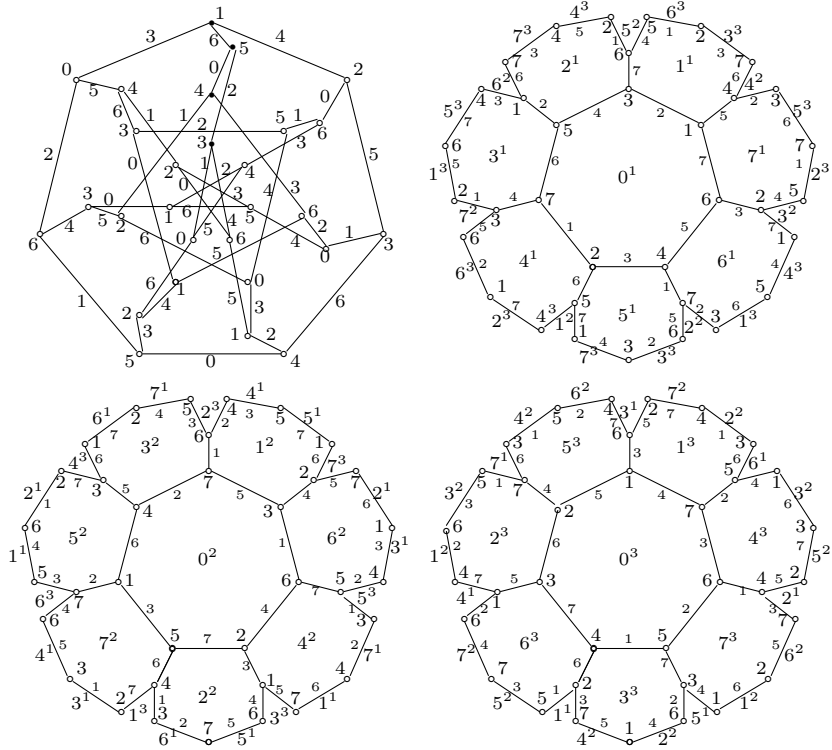
\begin{figure}[htp]
\unitlength=0.40mm
\special{em:linewidth 0.4pt}
\linethickness{0.4pt}
\begin{picture}(286.00,243.00)
\put(209.67,159.00){\circle{2.00}}
\put(233.67,159.00){\circle{2.00}}
\put(251.67,181.00){\circle{2.00}}
\put(191.67,181.00){\circle{2.00}}
\put(245.67,206.00){\circle{2.00}}
\put(197.67,206.00){\circle{2.00}}
\put(221.67,218.00){\circle{2.00}}
\emline{210.67}{159.00}{1}{232.67}{159.00}{2}
\emline{234.67}{160.00}{3}{250.67}{180.00}{4}
\emline{244.67}{207.00}{5}{222.67}{218.00}{6}
\emline{220.67}{218.00}{7}{198.67}{207.00}{8}
\emline{192.67}{180.00}{9}{208.67}{160.00}{10}
\put(239.67,147.00){\circle{2.00}}
\put(203.67,147.00){\circle{2.00}}
\put(177.67,178.00){\circle{2.00}}
\put(265.67,178.00){\circle{2.00}}
\put(256.67,215.00){\circle{2.00}}
\put(186.67,215.00){\circle{2.00}}
\put(221.67,230.00){\circle{2.00}}
\emline{221.67}{229.00}{11}{221.67}{219.00}{12}
\emline{255.67}{214.00}{13}{246.67}{207.00}{14}
\emline{187.67}{214.00}{15}{196.67}{207.00}{16}
\emline{190.67}{181.00}{17}{178.67}{178.00}{18}
\emline{252.67}{181.00}{19}{264.67}{178.00}{20}
\emline{233.67}{158.00}{21}{238.67}{148.00}{22}
\emline{209.67}{158.00}{23}{204.67}{148.00}{24}
\put(192.67,138.00){\circle{2.00}}
\emline{202.67}{146.00}{25}{193.67}{139.00}{26}
\put(166.67,169.00){\circle{2.00}}
\emline{176.67}{177.00}{27}{167.67}{170.00}{28}
\put(276.67,169.00){\circle{2.00}}
\emline{266.67}{177.00}{29}{275.67}{170.00}{30}
\put(250.67,138.00){\circle{2.00}}
\emline{240.67}{146.00}{31}{249.67}{139.00}{32}
\put(239.67,135.00){\circle{2.00}}
\emline{239.67}{146.00}{33}{239.67}{136.00}{34}
\put(203.67,135.00){\circle{2.00}}
\emline{203.67}{146.00}{35}{203.67}{136.00}{36}
\put(180.67,227.00){\circle{2.00}}
\emline{180.67}{226.00}{37}{185.67}{216.00}{38}
\put(279.67,181.00){\circle{2.00}}
\emline{278.67}{181.00}{39}{266.67}{178.00}{40}
\put(163.67,181.00){\circle{2.00}}
\emline{164.67}{181.00}{41}{176.67}{178.00}{42}
\put(262.67,227.00){\circle{2.00}}
\emline{262.67}{226.00}{43}{257.67}{216.00}{44}
\put(270.67,218.00){\circle{2.00}}
\emline{269.67}{218.00}{45}{257.67}{215.00}{46}
\put(172.67,218.00){\circle{2.00}}
\emline{173.67}{218.00}{47}{185.67}{215.00}{48}
\put(209.67,159.00){\circle{2.00}}
\put(227.67,242.00){\circle{2.00}}
\emline{227.67}{241.00}{49}{222.67}{231.00}{50}
\put(215.67,242.00){\circle{2.00}}
\emline{215.67}{241.00}{51}{220.67}{231.00}{52}
\put(221.67,128.00){\circle{2.00}}
\put(267.67,149.00){\circle{2.00}}
\put(175.67,149.00){\circle{2.00}}
\emline{166.67}{168.00}{53}{174.67}{150.00}{54}
\emline{176.67}{148.00}{55}{191.67}{138.00}{56}
\emline{204.67}{134.00}{57}{220.67}{128.00}{58}
\emline{222.67}{128.00}{59}{238.67}{134.00}{60}
\emline{251.67}{139.00}{61}{266.67}{148.00}{62}
\emline{268.67}{150.00}{63}{276.67}{168.00}{64}
\put(282.67,199.00){\circle{2.00}}
\put(160.67,199.00){\circle{2.00}}
\emline{160.67}{200.00}{65}{171.67}{217.00}{66}
\emline{271.67}{217.00}{67}{282.67}{200.00}{68}
\put(245.67,239.00){\circle{2.00}}
\put(197.67,239.00){\circle{2.00}}
\emline{181.67}{228.00}{69}{196.67}{239.00}{70}
\emline{198.67}{239.00}{71}{214.67}{242.00}{72}
\emline{228.67}{242.00}{73}{244.67}{239.00}{74}
\emline{246.67}{239.00}{75}{261.67}{228.00}{76}
\put(202.67,225.00){\makebox(0,0)[cc]{$^{2^1}$}}
\put(240.67,225.00){\makebox(0,0)[cc]{$^{1^1}$}}
\put(264.67,195.00){\makebox(0,0)[cc]{$^{7^1}$}}
\put(178.67,195.00){\makebox(0,0)[cc]{$^{3^1}$}}
\put(188.67,160.00){\makebox(0,0)[cc]{$^{4^1}$}}
\put(254.67,160.00){\makebox(0,0)[cc]{$^{6^1}$}}
\put(221.67,142.00){\makebox(0,0)[cc]{$^{5^1}$}}
\put(201.67,203.00){\makebox(0,0)[cc]{$^5$}}
\put(195.67,180.00){\makebox(0,0)[cc]{$^7$}}
\put(210.67,162.00){\makebox(0,0)[cc]{$^2$}}
\put(232.67,162.00){\makebox(0,0)[cc]{$^4$}}
\put(247.67,180.00){\makebox(0,0)[cc]{$^6$}}
\put(241.67,203.00){\makebox(0,0)[cc]{$^1$}}
\put(221.67,213.00){\makebox(0,0)[cc]{$^3$}}
\put(197.67,193.00){\makebox(0,0)[cc]{$^{_6}$}}
\put(202.67,171.00){\makebox(0,0)[cc]{$^{_1}$}}
\put(240.67,171.00){\makebox(0,0)[cc]{$^{_5}$}}
\put(221.67,161.00){\makebox(0,0)[cc]{$^{_3}$}}
\put(246.67,193.00){\makebox(0,0)[cc]{$^{_7}$}}
\put(232.67,209.00){\makebox(0,0)[cc]{$^{_2}$}}
\put(211.67,209.00){\makebox(0,0)[cc]{$^{_4}$}}
\put(253.67,208.00){\makebox(0,0)[cc]{$^{_5}$}}
\put(254.67,217.00){\makebox(0,0)[cc]{$^4$}}
\put(263.67,213.00){\makebox(0,0)[cc]{$^{_2}$}}
\put(270.67,213.00){\makebox(0,0)[cc]{$^3$}}
\put(274.67,206.00){\makebox(0,0)[cc]{$^{_6}$}}
\put(278.67,198.00){\makebox(0,0)[cc]{$^7$}}
\put(278.67,192.00){\makebox(0,0)[cc]{$^{_1}$}}
\put(276.67,183.00){\makebox(0,0)[cc]{$^5$}}
\put(271.67,181.00){\makebox(0,0)[cc]{$^{_4}$}}
\put(265.67,181.00){\makebox(0,0)[cc]{$^2$}}
\put(269.67,171.00){\makebox(0,0)[cc]{$^{_7}$}}
\put(272.67,167.00){\makebox(0,0)[cc]{$^1$}}
\put(270.67,159.00){\makebox(0,0)[cc]{$^{_4}$}}
\put(265.67,151.00){\makebox(0,0)[cc]{$^5$}}
\put(257.67,145.00){\makebox(0,0)[cc]{$^{_6}$}}
\put(250.67,141.00){\makebox(0,0)[cc]{$^3$}}
\put(245.67,144.00){\makebox(0,0)[cc]{$^{_2}$}}
\put(242.67,147.00){\makebox(0,0)[cc]{$^7$}}
\put(236.67,141.00){\makebox(0,0)[cc]{$^{_5}$}}
\put(236.67,136.00){\makebox(0,0)[cc]{$^6$}}
\put(229.67,133.00){\makebox(0,0)[cc]{$^{_2}$}}
\put(221.67,131.00){\makebox(0,0)[cc]{$^3$}}
\put(213.67,133.00){\makebox(0,0)[cc]{$^{_4}$}}
\put(206.67,136.00){\makebox(0,0)[cc]{$^1$}}
\put(206.67,141.00){\makebox(0,0)[cc]{$^{_7}$}}
\put(207.67,146.00){\makebox(0,0)[cc]{$^5$}}
\put(234.67,151.00){\makebox(0,0)[cc]{$^{_1}$}}
\put(204.67,153.00){\makebox(0,0)[cc]{$^{_6}$}}
\put(196.67,144.00){\makebox(0,0)[cc]{$^{_3}$}}
\put(192.67,141.00){\makebox(0,0)[cc]{$^4$}}
\put(184.67,145.00){\makebox(0,0)[cc]{$^{_7}$}}
\put(179.67,151.00){\makebox(0,0)[cc]{$^1$}}
\put(173.67,159.00){\makebox(0,0)[cc]{$^{_2}$}}
\put(173.67,171.00){\makebox(0,0)[cc]{$^{_5}$}}
\put(170.67,167.00){\makebox(0,0)[cc]{$^6$}}
\put(177.67,173.00){\makebox(0,0)[cc]{$^3$}}
\put(183.67,181.00){\makebox(0,0)[cc]{$^{_4}$}}
\put(171.67,181.00){\makebox(0,0)[cc]{$^{_1}$}}
\put(166.67,183.00){\makebox(0,0)[cc]{$^2$}}
\put(164.67,192.00){\makebox(0,0)[cc]{$^{_5}$}}
\put(164.67,198.00){\makebox(0,0)[cc]{$^6$}}
\put(168.67,206.00){\makebox(0,0)[cc]{$^{_7}$}}
\put(172.67,213.00){\makebox(0,0)[cc]{$^4$}}
\put(178.67,213.00){\makebox(0,0)[cc]{$^{_3}$}}
\put(184.67,210.00){\makebox(0,0)[cc]{$^1$}}
\put(193.67,212.00){\makebox(0,0)[cc]{$^{_2}$}}
\put(185.67,221.00){\makebox(0,0)[cc]{$^{_6}$}}
\put(184.67,226.00){\makebox(0,0)[cc]{$^7$}}
\put(190.67,229.00){\makebox(0,0)[cc]{$^{_3}$}}
\put(197.67,234.00){\makebox(0,0)[cc]{$^4$}}
\put(205.67,236.00){\makebox(0,0)[cc]{$^{_5}$}}
\put(214.67,237.00){\makebox(0,0)[cc]{$^2$}}
\put(217.67,233.00){\makebox(0,0)[cc]{$^{_1}$}}
\put(227.67,233.00){\makebox(0,0)[cc]{$^{_4}$}}
\put(229.67,237.00){\makebox(0,0)[cc]{$^5$}}
\put(235.67,236.00){\makebox(0,0)[cc]{$^{_1}$}}
\put(245.67,234.00){\makebox(0,0)[cc]{$^2$}}
\put(252.67,229.00){\makebox(0,0)[cc]{$^{_3}$}}
\put(258.67,226.00){\makebox(0,0)[cc]{$^7$}}
\put(257.67,221.00){\makebox(0,0)[cc]{$^{_6}$}}
\put(224.67,223.00){\makebox(0,0)[cc]{$^{_7}$}}
\put(221.67,188.00){\makebox(0,0)[cc]{$^{0^1}$}}
\put(179.67,219.00){\makebox(0,0)[cc]{$^{6^2}$}}
\put(222.67,237.00){\makebox(0,0)[cc]{$^{5^2}$}}
\put(263.67,219.00){\makebox(0,0)[cc]{$^{4^2}$}}
\put(274.67,175.00){\makebox(0,0)[cc]{$^{3^2}$}}
\put(243.67,138.00){\makebox(0,0)[cc]{$^{2^2}$}}
\put(200.67,139.00){\makebox(0,0)[cc]{$^{1^2}$}}
\put(168.67,175.00){\makebox(0,0)[cc]{$^{7^2}$}}
\put(237.67,243.00){\makebox(0,0)[cc]{$^{6^3}$}}
\put(205.67,243.00){\makebox(0,0)[cc]{$^{4^3}$}}
\put(186.67,235.00){\makebox(0,0)[cc]{$^{7^3}$}}
\put(256.67,235.00){\makebox(0,0)[cc]{$^{3^3}$}}
\put(279.67,210.00){\makebox(0,0)[cc]{$^{5^3}$}}
\put(285.67,189.00){\makebox(0,0)[cc]{$^{2^3}$}}
\put(276.67,156.00){\makebox(0,0)[cc]{$^{4^3}$}}
\put(261.67,140.00){\makebox(0,0)[cc]{$^{1^3}$}}
\put(234.67,128.00){\makebox(0,0)[cc]{$^{3^3}$}}
\put(208.67,128.00){\makebox(0,0)[cc]{$^{7^3}$}}
\put(179.67,140.00){\makebox(0,0)[cc]{$^{2^3}$}}
\put(167.67,156.00){\makebox(0,0)[cc]{$^{6^3}$}}
\put(158.67,189.00){\makebox(0,0)[cc]{$^{1^3}$}}
\put(163.67,210.00){\makebox(0,0)[cc]{$^{5^3}$}}
\put(70.00,32.00){\circle{2.00}}
\put(94.00,32.00){\circle{2.00}}
\put(112.00,54.00){\circle{2.00}}
\put(52.00,54.00){\circle{2.00}}
\put(106.00,79.00){\circle{2.00}}
\put(58.00,79.00){\circle{2.00}}
\put(82.00,91.00){\circle{2.00}}
\emline{71.00}{32.00}{77}{93.00}{32.00}{78}
\emline{95.00}{33.00}{79}{111.00}{53.00}{80}
\emline{105.00}{80.00}{81}{83.00}{91.00}{82}
\emline{81.00}{91.00}{83}{59.00}{80.00}{84}
\emline{53.00}{53.00}{85}{69.00}{33.00}{86}
\put(100.00,20.00){\circle{2.00}}
\put(64.00,20.00){\circle{2.00}}
\put(38.00,51.00){\circle{2.00}}
\put(126.00,51.00){\circle{2.00}}
\put(117.00,88.00){\circle{2.00}}
\put(47.00,88.00){\circle{2.00}}
\put(82.00,103.00){\circle{2.00}}
\emline{82.00}{102.00}{87}{82.00}{92.00}{88}
\emline{116.00}{87.00}{89}{107.00}{80.00}{90}
\emline{48.00}{87.00}{91}{57.00}{80.00}{92}
\emline{51.00}{54.00}{93}{39.00}{51.00}{94}
\emline{113.00}{54.00}{95}{125.00}{51.00}{96}
\emline{94.00}{31.00}{97}{99.00}{21.00}{98}
\emline{70.00}{31.00}{99}{65.00}{21.00}{100}
\put(53.00,11.00){\circle{2.00}}
\emline{63.00}{19.00}{101}{54.00}{12.00}{102}
\put(27.00,42.00){\circle{2.00}}
\emline{37.00}{50.00}{103}{28.00}{43.00}{104}
\put(137.00,42.00){\circle{2.00}}
\emline{127.00}{49.00}{105}{136.00}{42.00}{106}
\put(111.00,11.00){\circle{2.00}}
\emline{101.00}{19.00}{107}{110.00}{12.00}{108}
\put(100.00,8.00){\circle{2.00}}
\emline{100.00}{19.00}{109}{100.00}{9.00}{110}
\put(64.00,8.00){\circle{2.00}}
\emline{64.00}{19.00}{111}{64.00}{9.00}{112}
\put(41.00,100.00){\circle{2.00}}
\emline{41.00}{99.00}{113}{46.00}{89.00}{114}
\put(140.00,54.00){\circle{2.00}}
\emline{139.00}{54.00}{115}{127.00}{51.00}{116}
\put(24.00,54.00){\circle{2.00}}
\emline{25.00}{54.00}{117}{37.00}{51.00}{118}
\put(123.00,100.00){\circle{2.00}}
\emline{123.00}{99.00}{119}{118.00}{89.00}{120}
\put(131.00,91.00){\circle{2.00}}
\emline{130.00}{91.00}{121}{118.00}{88.00}{122}
\put(33.00,91.00){\circle{2.00}}
\emline{34.00}{91.00}{123}{46.00}{88.00}{124}
\put(70.00,32.00){\circle{2.00}}
\emline{70.00}{31.00}{125}{65.00}{21.00}{126}
\put(88.00,115.00){\circle{2.00}}
\emline{88.00}{114.00}{127}{83.00}{104.00}{128}
\put(76.00,115.00){\circle{2.00}}
\emline{76.00}{114.00}{129}{81.00}{104.00}{130}
\put(82.00,1.00){\circle{2.00}}
\put(128.00,22.00){\circle{2.00}}
\put(36.00,22.00){\circle{2.00}}
\emline{27.00}{41.00}{131}{35.00}{23.00}{132}
\emline{37.00}{21.00}{133}{52.00}{11.00}{134}
\emline{65.00}{7.00}{135}{81.00}{1.00}{136}
\emline{83.00}{1.00}{137}{99.00}{7.00}{138}
\emline{112.00}{12.00}{139}{127.00}{21.00}{140}
\emline{129.00}{23.00}{141}{137.00}{41.00}{142}
\put(143.00,72.00){\circle{2.00}}
\put(21.00,72.00){\circle{2.00}}
\emline{21.00}{73.00}{143}{32.00}{90.00}{144}
\emline{132.00}{90.00}{145}{143.00}{73.00}{146}
\put(106.00,112.00){\circle{2.00}}
\put(58.00,112.00){\circle{2.00}}
\emline{42.00}{101.00}{147}{57.00}{112.00}{148}
\emline{59.00}{112.00}{149}{75.00}{115.00}{150}
\emline{89.00}{115.00}{151}{105.00}{112.00}{152}
\emline{107.00}{112.00}{153}{122.00}{101.00}{154}
\put(82.00,1.00){\circle{2.00}}
\put(63.00,98.00){\makebox(0,0)[cc]{$^{3^2}$}}
\put(101.00,98.00){\makebox(0,0)[cc]{$^{1^2}$}}
\put(125.00,68.00){\makebox(0,0)[cc]{$^{6^2}$}}
\put(39.00,68.00){\makebox(0,0)[cc]{$^{5^2}$}}
\put(49.00,33.00){\makebox(0,0)[cc]{$^{7^2}$}}
\put(115.00,33.00){\makebox(0,0)[cc]{$^{4^2}$}}
\put(82.00,15.00){\makebox(0,0)[cc]{$^{2^2}$}}
\put(62.00,76.00){\makebox(0,0)[cc]{$^4$}}
\put(56.00,53.00){\makebox(0,0)[cc]{$^1$}}
\put(71.00,35.00){\makebox(0,0)[cc]{$^5$}}
\put(93.00,35.00){\makebox(0,0)[cc]{$^2$}}
\put(108.00,53.00){\makebox(0,0)[cc]{$^6$}}
\put(102.00,76.00){\makebox(0,0)[cc]{$^3$}}
\put(82.00,86.00){\makebox(0,0)[cc]{$^7$}}
\put(58.00,66.00){\makebox(0,0)[cc]{$^{_6}$}}
\put(63.00,44.00){\makebox(0,0)[cc]{$^{_3}$}}
\put(101.00,44.00){\makebox(0,0)[cc]{$^{_4}$}}
\put(82.00,34.00){\makebox(0,0)[cc]{$^{_7}$}}
\put(107.00,66.00){\makebox(0,0)[cc]{$^{_1}$}}
\put(93.00,82.00){\makebox(0,0)[cc]{$^{_5}$}}
\put(72.00,82.00){\makebox(0,0)[cc]{$^{_2}$}}
\put(115.00,82.00){\makebox(0,0)[cc]{$^{_4}$}}
\put(114.00,89.00){\makebox(0,0)[cc]{$^2$}}
\put(124.00,86.00){\makebox(0,0)[cc]{$^{_5}$}}
\put(131.00,86.00){\makebox(0,0)[cc]{$^7$}}
\put(135.00,79.00){\makebox(0,0)[cc]{$^{_6}$}}
\put(139.00,71.00){\makebox(0,0)[cc]{$^1$}}
\put(139.00,65.00){\makebox(0,0)[cc]{$^{_3}$}}
\put(137.00,56.00){\makebox(0,0)[cc]{$^4$}}
\put(132.00,54.00){\makebox(0,0)[cc]{$^{_2}$}}
\put(126.00,54.00){\makebox(0,0)[cc]{$^5$}}
\put(130.00,44.00){\makebox(0,0)[cc]{$^{_1}$}}
\put(133.00,40.00){\makebox(0,0)[cc]{$^3$}}
\put(131.00,32.00){\makebox(0,0)[cc]{$^{_2}$}}
\put(126.00,24.00){\makebox(0,0)[cc]{$^4$}}
\put(118.00,18.00){\makebox(0,0)[cc]{$^{_6}$}}
\put(111.00,14.00){\makebox(0,0)[cc]{$^7$}}
\put(106.00,17.00){\makebox(0,0)[cc]{$^{_5}$}}
\put(103.00,20.00){\makebox(0,0)[cc]{$^1$}}
\put(97.00,14.00){\makebox(0,0)[cc]{$^{_4}$}}
\put(97.00,9.00){\makebox(0,0)[cc]{$^6$}}
\put(90.00,6.00){\makebox(0,0)[cc]{$^{_5}$}}
\put(82.00,4.00){\makebox(0,0)[cc]{$^7$}}
\put(74.00,6.00){\makebox(0,0)[cc]{$^{_2}$}}
\put(67.00,9.00){\makebox(0,0)[cc]{$^3$}}
\put(67.00,14.00){\makebox(0,0)[cc]{$^{_1}$}}
\put(68.00,19.00){\makebox(0,0)[cc]{$^4$}}
\put(95.00,25.00){\makebox(0,0)[cc]{$^{_3}$}}
\put(65.00,26.00){\makebox(0,0)[cc]{$^{_6}$}}
\put(57.00,17.00){\makebox(0,0)[cc]{$^{_7}$}}
\put(53.00,14.00){\makebox(0,0)[cc]{$^2$}}
\put(45.00,18.00){\makebox(0,0)[cc]{$^{_1}$}}
\put(40.00,24.00){\makebox(0,0)[cc]{$^3$}}
\put(34.00,32.00){\makebox(0,0)[cc]{$^{_5}$}}
\put(35.00,44.00){\makebox(0,0)[cc]{$^{_4}$}}
\put(31.00,40.00){\makebox(0,0)[cc]{$^6$}}
\put(39.00,46.00){\makebox(0,0)[cc]{$^7$}}
\put(44.00,54.00){\makebox(0,0)[cc]{$^{_2}$}}
\put(32.00,54.00){\makebox(0,0)[cc]{$^{_3}$}}
\put(27.00,56.00){\makebox(0,0)[cc]{$^5$}}
\put(25.00,65.00){\makebox(0,0)[cc]{$^{_4}$}}
\put(25.00,71.00){\makebox(0,0)[cc]{$^6$}}
\put(29.00,79.00){\makebox(0,0)[cc]{$^{_1}$}}
\put(33.00,86.00){\makebox(0,0)[cc]{$^2$}}
\put(39.00,86.00){\makebox(0,0)[cc]{$^{_7}$}}
\put(45.00,83.00){\makebox(0,0)[cc]{$^3$}}
\put(55.00,84.00){\makebox(0,0)[cc]{$^{_5}$}}
\put(46.00,94.00){\makebox(0,0)[cc]{$^{_6}$}}
\put(45.00,99.00){\makebox(0,0)[cc]{$^1$}}
\put(51.00,102.00){\makebox(0,0)[cc]{$^{_7}$}}
\put(58.00,107.00){\makebox(0,0)[cc]{$^2$}}
\put(66.00,109.00){\makebox(0,0)[cc]{$^{_4}$}}
\put(75.00,110.00){\makebox(0,0)[cc]{$^5$}}
\put(77.00,106.00){\makebox(0,0)[cc]{$^{_3}$}}
\put(79.00,101.00){\makebox(0,0)[cc]{$^6$}}
\put(88.00,106.00){\makebox(0,0)[cc]{$^{_2}$}}
\put(90.00,110.00){\makebox(0,0)[cc]{$^4$}}
\put(96.00,109.00){\makebox(0,0)[cc]{$^{_3}$}}
\put(106.00,107.00){\makebox(0,0)[cc]{$^5$}}
\put(113.00,102.00){\makebox(0,0)[cc]{$^{_7}$}}
\put(119.00,99.00){\makebox(0,0)[cc]{$^1$}}
\put(118.00,94.00){\makebox(0,0)[cc]{$^{_6}$}}
\put(85.00,96.00){\makebox(0,0)[cc]{$^{_1}$}}
\put(82.00,61.00){\makebox(0,0)[cc]{$^{0^2}$}}
\put(39.00,93.00){\makebox(0,0)[cc]{$^{4^3}$}}
\put(83.00,110.00){\makebox(0,0)[cc]{$^{2^3}$}}
\put(125.00,92.00){\makebox(0,0)[cc]{$^{7^3}$}}
\put(135.00,48.00){\makebox(0,0)[cc]{$^{5^3}$}}
\put(106.00,8.00){\makebox(0,0)[cc]{$^{3^3}$}}
\put(59.00,9.00){\makebox(0,0)[cc]{$^{1^3}$}}
\put(28.00,47.00){\makebox(0,0)[cc]{$^{6^3}$}}
\put(98.00,116.00){\makebox(0,0)[cc]{$^{4^1}$}}
\put(66.00,116.00){\makebox(0,0)[cc]{$^{7^1}$}}
\put(47.00,108.00){\makebox(0,0)[cc]{$^{6^1}$}}
\put(117.00,108.00){\makebox(0,0)[cc]{$^{5^1}$}}
\put(140.00,83.00){\makebox(0,0)[cc]{$^{2^1}$}}
\put(146.00,64.96){\makebox(0,0)[cc]{$^{3^1}$}}
\put(137.00,29.00){\makebox(0,0)[cc]{$^{7^1}$}}
\put(122.00,13.00){\makebox(0,0)[cc]{$^{1^1}$}}
\put(95.00,1.00){\makebox(0,0)[cc]{$^{5^1}$}}
\put(69.00,1.00){\makebox(0,0)[cc]{$^{6^1}$}}
\put(40.00,13.00){\makebox(0,0)[cc]{$^{3^1}$}}
\put(28.00,29.00){\makebox(0,0)[cc]{$^{4^1}$}}
\put(19.00,62.00){\makebox(0,0)[cc]{$^{1^1}$}}
\put(24.00,83.00){\makebox(0,0)[cc]{$^{2^1}$}}
\put(210.00,32.00){\circle{2.00}}
\put(234.00,32.00){\circle{2.00}}
\put(252.00,54.00){\circle{2.00}}
\put(192.00,54.00){\circle{2.00}}
\put(246.00,79.00){\circle{2.00}}
\put(198.00,78.00){\circle{2.00}}
\put(222.00,91.00){\circle{2.00}}
\emline{211.00}{32.00}{155}{233.00}{32.00}{156}
\emline{235.00}{33.00}{157}{251.00}{53.00}{158}
\emline{245.00}{80.00}{159}{223.00}{91.00}{160}
\emline{221.00}{91.00}{161}{199.00}{80.00}{162}
\emline{193.00}{53.00}{163}{209.00}{33.00}{164}
\put(240.00,20.00){\circle{2.00}}
\put(204.00,20.00){\circle{2.00}}
\put(178.00,51.00){\circle{2.00}}
\put(266.00,51.00){\circle{2.00}}
\put(257.00,88.00){\circle{2.00}}
\put(187.00,88.00){\circle{2.00}}
\put(222.00,103.00){\circle{2.00}}
\emline{222.00}{102.00}{165}{222.00}{92.00}{166}
\emline{256.00}{87.00}{167}{247.00}{80.00}{168}
\emline{188.00}{87.00}{169}{197.00}{80.00}{170}
\emline{191.00}{54.00}{171}{179.00}{51.00}{172}
\emline{253.00}{54.00}{173}{265.00}{51.00}{174}
\emline{234.00}{31.00}{175}{239.00}{21.00}{176}
\emline{210.00}{31.00}{177}{205.00}{21.00}{178}
\put(193.00,11.00){\circle{2.00}}
\emline{203.00}{18.00}{179}{194.00}{11.00}{180}
\put(167.00,42.00){\circle{2.00}}
\emline{177.00}{50.00}{181}{168.00}{43.00}{182}
\put(277.00,42.00){\circle{2.00}}
\emline{267.00}{49.00}{183}{276.00}{42.00}{184}
\put(251.00,11.00){\circle{2.00}}
\emline{241.00}{18.00}{185}{250.00}{11.00}{186}
\put(240.00,8.00){\circle{2.00}}
\emline{240.00}{19.00}{187}{240.00}{9.00}{188}
\put(204.00,8.00){\circle{2.00}}
\emline{204.00}{19.00}{189}{204.00}{9.00}{190}
\put(181.00,100.00){\circle{2.00}}
\emline{181.00}{99.00}{191}{186.00}{89.00}{192}
\put(280.00,54.00){\circle{2.00}}
\emline{279.00}{54.00}{193}{267.00}{51.00}{194}
\put(164.00,54.00){\circle{2.00}}
\emline{165.00}{54.00}{195}{177.00}{51.00}{196}
\put(263.00,100.00){\circle{2.00}}
\emline{263.00}{99.00}{197}{258.00}{89.00}{198}
\put(271.00,91.00){\circle{2.00}}
\emline{270.00}{91.00}{199}{258.00}{88.00}{200}
\put(173.00,91.00){\circle{2.00}}
\emline{174.00}{91.00}{201}{186.00}{88.00}{202}
\put(210.00,32.00){\circle{2.00}}
\emline{210.00}{31.00}{203}{205.00}{21.00}{204}
\put(228.00,115.00){\circle{2.00}}
\emline{228.00}{114.00}{205}{223.00}{104.00}{206}
\put(216.00,115.00){\circle{2.00}}
\emline{216.00}{114.00}{207}{221.00}{104.00}{208}
\put(222.00,1.00){\circle{2.00}}
\put(268.00,22.00){\circle{2.00}}
\put(176.00,22.00){\circle{2.00}}
\emline{167.00}{41.00}{209}{175.00}{23.00}{210}
\emline{177.00}{21.00}{211}{192.00}{11.00}{212}
\emline{205.00}{7.00}{213}{221.00}{1.00}{214}
\emline{223.00}{1.00}{215}{239.00}{7.00}{216}
\emline{252.00}{12.00}{217}{267.00}{21.00}{218}
\emline{269.00}{23.00}{219}{277.00}{41.00}{220}
\put(283.00,72.00){\circle{2.00}}
\put(161.00,71.00){\circle{2.00}}
\emline{161.00}{73.00}{221}{172.00}{90.00}{222}
\emline{272.00}{90.00}{223}{283.00}{73.00}{224}
\put(246.00,112.00){\circle{2.00}}
\put(198.00,112.00){\circle{2.00}}
\emline{182.00}{101.00}{225}{197.00}{112.00}{226}
\emline{199.00}{112.00}{227}{215.00}{115.00}{228}
\emline{229.00}{115.00}{229}{245.00}{112.00}{230}
\emline{247.00}{112.00}{231}{262.00}{101.00}{232}
\put(222.00,1.00){\circle{2.00}}
\put(203.00,98.00){\makebox(0,0)[cc]{$^{5^3}$}}
\put(241.00,98.00){\makebox(0,0)[cc]{$^{1^3}$}}
\put(265.00,68.00){\makebox(0,0)[cc]{$^{4^3}$}}
\put(179.00,68.00){\makebox(0,0)[cc]{$^{2^3}$}}
\put(189.00,33.00){\makebox(0,0)[cc]{$^{6^3}$}}
\put(255.00,33.00){\makebox(0,0)[cc]{$^{7^3}$}}
\put(222.00,15.00){\makebox(0,0)[cc]{$^{3^3}$}}
\put(202.00,76.00){\makebox(0,0)[cc]{$^2$}}
\put(196.00,53.00){\makebox(0,0)[cc]{$^3$}}
\put(211.00,35.00){\makebox(0,0)[cc]{$^4$}}
\put(233.00,35.00){\makebox(0,0)[cc]{$^5$}}
\put(248.00,53.00){\makebox(0,0)[cc]{$^6$}}
\put(242.00,76.00){\makebox(0,0)[cc]{$^7$}}
\put(222.00,86.00){\makebox(0,0)[cc]{$^1$}}
\put(198.00,66.00){\makebox(0,0)[cc]{$^{_6}$}}
\put(203.00,44.00){\makebox(0,0)[cc]{$^{_7}$}}
\put(241.00,44.00){\makebox(0,0)[cc]{$^{_2}$}}
\put(222.00,34.00){\makebox(0,0)[cc]{$^{_1}$}}
\put(247.00,66.00){\makebox(0,0)[cc]{$^{_3}$}}
\put(233.00,82.00){\makebox(0,0)[cc]{$^{_4}$}}
\put(212.00,82.00){\makebox(0,0)[cc]{$^{_5}$}}
\put(254.00,81.00){\makebox(0,0)[cc]{$^{_2}$}}
\put(255.00,90.00){\makebox(0,0)[cc]{$^5$}}
\put(264.00,86.00){\makebox(0,0)[cc]{$^{_4}$}}
\put(271.00,86.00){\makebox(0,0)[cc]{$^1$}}
\put(275.00,79.00){\makebox(0,0)[cc]{$^{_6}$}}
\put(279.00,71.00){\makebox(0,0)[cc]{$^3$}}
\put(279.00,65.00){\makebox(0,0)[cc]{$^{_7}$}}
\put(277.00,56.00){\makebox(0,0)[cc]{$^2$}}
\put(272.00,54.00){\makebox(0,0)[cc]{$^{_5}$}}
\put(266.00,54.00){\makebox(0,0)[cc]{$^4$}}
\put(260.00,50.00){\makebox(0,0)[cc]{$^{_1}$}}
\put(270.00,44.00){\makebox(0,0)[cc]{$^{_3}$}}
\put(273.00,40.00){\makebox(0,0)[cc]{$^7$}}
\put(271.00,32.00){\makebox(0,0)[cc]{$^{_5}$}}
\put(266.00,24.00){\makebox(0,0)[cc]{$^2$}}
\put(258.00,18.00){\makebox(0,0)[cc]{$^{_6}$}}
\put(251.00,14.00){\makebox(0,0)[cc]{$^1$}}
\put(246.00,17.00){\makebox(0,0)[cc]{$^{_4}$}}
\put(243.00,20.00){\makebox(0,0)[cc]{$^3$}}
\put(237.00,14.00){\makebox(0,0)[cc]{$^{_2}$}}
\put(237.00,9.00){\makebox(0,0)[cc]{$^6$}}
\put(230.00,6.00){\makebox(0,0)[cc]{$^{_4}$}}
\put(222.00,4.00){\makebox(0,0)[cc]{$^1$}}
\put(214.00,6.00){\makebox(0,0)[cc]{$^{_5}$}}
\put(207.00,9.00){\makebox(0,0)[cc]{$^7$}}
\put(207.00,14.00){\makebox(0,0)[cc]{$^{_3}$}}
\put(208.00,19.00){\makebox(0,0)[cc]{$^2$}}
\put(234.00,26.00){\makebox(0,0)[cc]{$^{_7}$}}
\put(205.00,26.00){\makebox(0,0)[cc]{$^{_6}$}}
\put(198.00,17.00){\makebox(0,0)[cc]{$^{_1}$}}
\put(193.00,14.00){\makebox(0,0)[cc]{$^5$}}
\put(185.00,18.00){\makebox(0,0)[cc]{$^{_3}$}}
\put(180.00,24.00){\makebox(0,0)[cc]{$^7$}}
\put(174.00,32.00){\makebox(0,0)[cc]{$^{_4}$}}
\put(174.00,44.00){\makebox(0,0)[cc]{$^{_2}$}}
\put(171.00,40.00){\makebox(0,0)[cc]{$^6$}}
\put(178.00,46.00){\makebox(0,0)[cc]{$^1$}}
\put(184.00,54.00){\makebox(0,0)[cc]{$^{_5}$}}
\put(172.00,54.00){\makebox(0,0)[cc]{$^{_7}$}}
\put(167.00,56.00){\makebox(0,0)[cc]{$^4$}}
\put(165.00,65.00){\makebox(0,0)[cc]{$^{_2}$}}
\put(165.00,71.00){\makebox(0,0)[cc]{$^6$}}
\put(169.00,79.00){\makebox(0,0)[cc]{$^{_3}$}}
\put(173.00,86.00){\makebox(0,0)[cc]{$^5$}}
\put(179.00,86.00){\makebox(0,0)[cc]{$^{_1}$}}
\put(185.00,83.00){\makebox(0,0)[cc]{$^7$}}
\put(195.00,84.00){\makebox(0,0)[cc]{$^{_4}$}}
\put(186.00,94.00){\makebox(0,0)[cc]{$^{_6}$}}
\put(185.00,99.00){\makebox(0,0)[cc]{$^3$}}
\put(191.00,102.00){\makebox(0,0)[cc]{$^{_1}$}}
\put(198.00,107.00){\makebox(0,0)[cc]{$^5$}}
\put(206.00,109.00){\makebox(0,0)[cc]{$^{_2}$}}
\put(215.00,110.00){\makebox(0,0)[cc]{$^4$}}
\put(217.00,106.00){\makebox(0,0)[cc]{$^{_7}$}}
\put(218.00,102.00){\makebox(0,0)[cc]{$^6$}}
\put(228.00,106.00){\makebox(0,0)[cc]{$^{_5}$}}
\put(230.00,110.00){\makebox(0,0)[cc]{$^2$}}
\put(236.00,109.00){\makebox(0,0)[cc]{$^{_7}$}}
\put(246.00,107.00){\makebox(0,0)[cc]{$^4$}}
\put(253.00,102.00){\makebox(0,0)[cc]{$^{_1}$}}
\put(259.00,99.00){\makebox(0,0)[cc]{$^3$}}
\put(258.00,94.00){\makebox(0,0)[cc]{$^{_6}$}}
\put(225.00,96.00){\makebox(0,0)[cc]{$^{_3}$}}
\put(222.00,61.00){\makebox(0,0)[cc]{$^{0^3}$}}
\put(179.00,93.00){\makebox(0,0)[cc]{$^{7^1}$}}
\put(222.00,110.00){\makebox(0,0)[cc]{$^{3^1}$}}
\put(265.00,93.00){\makebox(0,0)[cc]{$^{6^1}$}}
\put(276.00,48.00){\makebox(0,0)[cc]{$^{2^1}$}}
\put(246.00,8.00){\makebox(0,0)[cc]{$^{5^1}$}}
\put(200.00,9.00){\makebox(0,0)[cc]{$^{1^1}$}}
\put(168.00,48.00){\makebox(0,0)[cc]{$^{4^1}$}}
\put(238.00,116.00){\makebox(0,0)[cc]{$^{7^2}$}}
\put(206.00,116.00){\makebox(0,0)[cc]{$^{6^2}$}}
\put(187.00,108.00){\makebox(0,0)[cc]{$^{4^2}$}}
\put(257.00,108.00){\makebox(0,0)[cc]{$^{2^2}$}}
\put(280.00,83.00){\makebox(0,0)[cc]{$^{3^2}$}}
\put(286.00,62.00){\makebox(0,0)[cc]{$^{5^2}$}}
\put(277.00,29.00){\makebox(0,0)[cc]{$^{6^2}$}}
\put(262.00,13.00){\makebox(0,0)[cc]{$^{1^2}$}}
\put(235.00,1.00){\makebox(0,0)[cc]{$^{2^2}$}}
\put(208.00,1.00){\makebox(0,0)[cc]{$^{4^2}$}}
\put(180.00,13.00){\makebox(0,0)[cc]{$^{5^2}$}}
\put(168.00,29.00){\makebox(0,0)[cc]{$^{7^2}$}}
\put(159.00,62.00){\makebox(0,0)[cc]{$^{1^2}$}}
\put(164.00,83.00){\makebox(0,0)[cc]{$^{3^2}$}}
\emline{24.00}{55.00}{233}{21.00}{71.00}{234}
\emline{58.00}{78.00}{235}{52.00}{55.00}{236}
\emline{106.00}{78.00}{237}{112.00}{55.00}{238}
\emline{143.00}{71.00}{239}{140.00}{55.00}{240}
\emline{164.00}{55.00}{241}{161.00}{71.00}{242}
\emline{198.00}{78.00}{243}{192.00}{55.00}{244}
\emline{246.00}{78.00}{245}{252.00}{55.00}{246}
\emline{280.00}{55.00}{247}{283.00}{71.00}{248}
\emline{161.00}{198.00}{249}{164.00}{182.00}{250}
\emline{280.00}{182.00}{251}{283.00}{198.00}{252}
\emline{192.00}{182.00}{253}{198.00}{205.00}{254}
\emline{252.00}{182.00}{255}{246.00}{205.00}{256}
\put(71.00,154.00){\circle{2.00}}
\put(95.00,154.00){\circle{2.00}}
\put(53.00,176.00){\circle{2.00}}
\put(107.00,204.00){\circle{2.00}}
\put(71.00,193.00){\circle{2.00}}
\put(89.00,168.00){\circle{2.00}}
\put(77.00,168.00){\circle{2.00}}
\put(42.00,179.00){\circle{2.00}}
\put(113.00,176.00){\circle{2.00}}
\put(94.00,193.00){\circle{2.00}}
\put(71.00,154.00){\circle{2.00}}
\put(38.00,221.00){\circle{2.00}}
\put(128.00,221.00){\circle{2.00}}
\put(26.00,170.00){\circle{2.00}}
\put(140.00,170.00){\circle{2.00}}
\put(59.00,130.00){\circle{2.00}}
\put(107.00,130.00){\circle{2.00}}
\emline{84.00}{240.00}{257}{127.00}{222.00}{258}
\emline{39.00}{222.00}{259}{82.00}{240.00}{260}
\emline{128.00}{220.00}{261}{140.00}{171.00}{262}
\emline{38.00}{220.00}{263}{26.00}{171.00}{264}
\emline{26.00}{169.00}{265}{58.00}{131.00}{266}
\emline{60.00}{130.00}{267}{106.00}{130.00}{268}
\emline{108.00}{131.00}{269}{140.00}{169.00}{270}
\put(58.00,204.00){\circle{2.00}}
\emline{84.00}{216.00}{271}{114.00}{177.00}{272}
\emline{113.00}{175.00}{273}{72.00}{154.00}{274}
\emline{70.00}{155.00}{275}{57.00}{203.00}{276}
\emline{59.00}{204.00}{277}{106.00}{204.00}{278}
\emline{108.00}{203.00}{279}{96.00}{154.00}{280}
\emline{94.00}{154.00}{281}{53.00}{175.00}{282}
\emline{52.00}{177.00}{283}{82.00}{216.00}{284}
\put(96.00,179.00){\circle{2.00}}
\emline{83.00}{199.00}{285}{89.00}{169.00}{286}
\emline{89.00}{169.00}{287}{72.00}{192.00}{288}
\emline{97.00}{178.00}{289}{120.00}{165.00}{290}
\emline{68.00}{179.00}{291}{43.00}{179.00}{292}
\emline{70.00}{179.00}{293}{93.00}{192.00}{294}
\emline{93.00}{192.00}{295}{77.00}{169.00}{296}
\emline{77.00}{169.00}{297}{83.00}{199.00}{298}
\put(119.00,207.00){\circle{2.00}}
\emline{95.00}{193.00}{299}{118.00}{206.00}{300}
\put(121.00,165.00){\circle{2.00}}
\emline{72.00}{192.00}{301}{95.00}{179.00}{302}
\emline{89.00}{167.00}{303}{95.00}{137.00}{304}
\put(95.00,136.00){\circle{2.00}}
\put(59.00,143.00){\circle{2.00}}
\emline{76.00}{167.00}{305}{60.00}{144.00}{306}
\put(69.00,179.00){\circle{2.00}}
\emline{95.00}{179.00}{307}{70.00}{179.00}{308}
\put(53.00,218.00){\circle{2.00}}
\emline{71.00}{194.00}{309}{54.00}{217.00}{310}
\emline{83.00}{239.00}{311}{89.00}{232.00}{312}
\emline{127.00}{220.00}{313}{120.00}{208.00}{314}
\emline{139.00}{170.00}{315}{122.00}{166.00}{316}
\emline{106.00}{131.00}{317}{96.00}{136.00}{318}
\emline{60.00}{131.00}{319}{58.00}{142.00}{320}
\emline{27.00}{170.00}{321}{41.00}{178.00}{322}
\emline{39.00}{221.00}{323}{52.00}{218.00}{324}
\emline{83.00}{201.00}{325}{91.00}{231.00}{326}
\emline{83.00}{217.00}{327}{89.00}{231.00}{328}
\emline{108.00}{205.00}{329}{118.00}{208.00}{330}
\emline{114.00}{175.00}{331}{121.00}{166.00}{332}
\emline{96.00}{153.00}{333}{96.00}{137.00}{334}
\emline{60.00}{142.00}{335}{71.00}{153.00}{336}
\emline{52.00}{217.00}{337}{57.00}{205.00}{338}
\put(86.00,242.00){\makebox(0,0)[cc]{$^1$}}
\put(106.00,234.00){\makebox(0,0)[cc]{$^4$}}
\put(132.00,223.00){\makebox(0,0)[cc]{$^2$}}
\put(142.00,165.00){\makebox(0,0)[cc]{$^3$}}
\put(111.00,126.00){\makebox(0,0)[cc]{$^4$}}
\put(126.00,146.00){\makebox(0,0)[cc]{$^6$}}
\put(137.00,195.00){\makebox(0,0)[cc]{$^5$}}
\put(83.00,126.00){\makebox(0,0)[cc]{$^0$}}
\put(40.00,146.00){\makebox(0,0)[cc]{$^1$}}
\put(24.00,165.00){\makebox(0,0)[cc]{$^6$}}
\put(29.00,195.00){\makebox(0,0)[cc]{$^2$}}
\put(35.00,223.00){\makebox(0,0)[cc]{$^0$}}
\put(63.00,234.00){\makebox(0,0)[cc]{$^3$}}
\put(84.00,233.00){\makebox(0,0)[cc]{$^6$}}
\put(94.00,230.00){\makebox(0,0)[cc]{$^5$}}
\put(79.00,216.00){\makebox(0,0)[cc]{$^4$}}
\put(84.00,224.00){\makebox(0,0)[cc]{$^0$}}
\put(90.00,216.00){\makebox(0,0)[cc]{$^2$}}
\put(80.00,199.00){\makebox(0,0)[cc]{$^3$}}
\put(96.00,189.00){\makebox(0,0)[cc]{$^4$}}
\put(96.00,174.00){\makebox(0,0)[cc]{$^5$}}
\put(93.00,167.00){\makebox(0,0)[cc]{$^6$}}
\put(73.00,167.00){\makebox(0,0)[cc]{$^0$}}
\put(69.00,174.00){\makebox(0,0)[cc]{$^1$}}
\put(69.00,189.00){\makebox(0,0)[cc]{$^2$}}
\put(81.00,192.00){\makebox(0,0)[cc]{$^1$}}
\put(74.00,185.00){\makebox(0,0)[cc]{$^0$}}
\put(75.00,176.00){\makebox(0,0)[cc]{$^6$}}
\put(82.00,170.00){\makebox(0,0)[cc]{$^5$}}
\put(90.00,173.00){\makebox(0,0)[cc]{$^4$}}
\put(92.00,182.00){\makebox(0,0)[cc]{$^3$}}
\put(88.00,190.00){\makebox(0,0)[cc]{$^2$}}
\put(112.00,208.00){\makebox(0,0)[cc]{$^1$}}
\put(106.00,207.00){\makebox(0,0)[cc]{$^5$}}
\put(120.00,202.00){\makebox(0,0)[cc]{$^6$}}
\put(121.00,214.00){\makebox(0,0)[cc]{$^0$}}
\put(63.00,208.00){\makebox(0,0)[cc]{$^1$}}
\put(80.00,205.00){\makebox(0,0)[cc]{$^2$}}
\put(57.00,217.00){\makebox(0,0)[cc]{$^4$}}
\put(52.00,209.00){\makebox(0,0)[cc]{$^6$}}
\put(54.00,203.00){\makebox(0,0)[cc]{$^3$}}
\put(57.00,192.00){\makebox(0,0)[cc]{$^0$}}
\put(73.00,150.00){\makebox(0,0)[cc]{$^1$}}
\put(55.00,142.00){\makebox(0,0)[cc]{$^2$}}
\put(63.00,152.00){\makebox(0,0)[cc]{$^6$}}
\put(62.00,136.00){\makebox(0,0)[cc]{$^3$}}
\put(91.00,135.00){\makebox(0,0)[cc]{$^1$}}
\put(103.00,135.00){\makebox(0,0)[cc]{$^2$}}
\put(121.00,161.00){\makebox(0,0)[cc]{$^0$}}
\put(129.00,169.00){\makebox(0,0)[cc]{$^1$}}
\put(120.00,171.00){\makebox(0,0)[cc]{$^2$}}
\put(117.00,176.00){\makebox(0,0)[cc]{$^6$}}
\put(111.00,184.00){\makebox(0,0)[cc]{$^3$}}
\put(53.00,171.00){\makebox(0,0)[cc]{$^2$}}
\put(61.00,167.00){\makebox(0,0)[cc]{$^6$}}
\put(47.00,173.00){\makebox(0,0)[cc]{$^5$}}
\put(42.00,182.00){\makebox(0,0)[cc]{$^3$}}
\put(49.00,180.00){\makebox(0,0)[cc]{$^0$}}
\put(111.00,166.00){\makebox(0,0)[cc]{$^4$}}
\put(74.00,208.00){\makebox(0,0)[cc]{$^1$}}
\put(68.00,145.00){\makebox(0,0)[cc]{$^4$}}
\put(99.00,144.00){\makebox(0,0)[cc]{$^3$}}
\put(99.00,153.00){\makebox(0,0)[cc]{$^0$}}
\put(85.00,162.00){\makebox(0,0)[cc]{$^5$}}
\put(44.00,216.00){\makebox(0,0)[cc]{$^5$}}
\put(36.00,171.00){\makebox(0,0)[cc]{$^4$}}
\put(113.00,199.00){\makebox(0,0)[cc]{$^3$}}
\put(90.00,147.00){\makebox(0,0)[cc]{$^5$}}
\put(102.00,183.00){\makebox(0,0)[cc]{$^4$}}
\put(56.00,126.00){\makebox(0,0)[cc]{$^5$}}
\emline{43.00}{179.00}{339}{52.00}{176.00}{340}
\put(83.00,240.00){\circle*{2.00}}
\put(90.00,232.00){\circle*{2.00}}
\put(83.00,216.00){\circle*{2.00}}
\put(83.00,200.00){\circle*{2.00}}
\put(258.67,176.00){\makebox(0,0)[cc]{$^{_3}$}}
\put(218.67,228.00){\makebox(0,0)[cc]{$^6$}}
\put(118.00,50.00){\makebox(0,0)[cc]{$^{_7}$}}
\end{picture}
\caption{$\mathcal F$-colored $\G$ and the three charts of $\G'$}
\end{figure}

\noindent Continuing this way, the oriented 7-cycles $\vec{C}_7^2$, indicated by means of the symbols $i^j$ (corresponding respectively to their square-root cycles $\vec{C}_7=\underline{i^j}$), where $i\in\{0\}\cup J_7$ and $j\in J_3=\{1,2,3\}$, are presented now as follows, by means of the colors $c(u_i)$ for their composing vertices and arcs, that make them pairwise distinguishable, as claimed, thus providing a distinctive and well-defined notation for them:
$$\begin{array}{cccccc}
^{0^1:}_{1^1:} & ^{(1_23_45_67_12_34_56_7);}_{(1_54_67_32_15_46_73_2);} & ^{0^2:}_{1^2:} & ^{(1_35_72_46_13_57_24_6);}_{(1_75_34_26_17_53_42_6);} & ^{0^3:}_{1^3:} & ^{(1_52_63_74_15_26_37_4);}_{(1_47_25_63_14_72_56_3);} \\
^{2^1:}_{3^1:} & ^{(1_25_43_76_12_54_37_6);}_{(1_34_76_52_13_47_65_2);} & ^{2^2:}_{3^2:} & ^{(1_32_75_64_13_27_56_4);}_{(1_63_54_27_16_35_42_7);} & ^{2^3:}_{3^3:} & ^{(1_53_62_47_15_36_24_7);}_{(1_46_23_75_14_62_37_5);} \\
^{4^1:}_{5^1:} & ^{(1_74_35_62_17_43_56_2);}_{(1_43_26_57_14_32_65_7);} & ^{4^2:}_{5^2:} & ^{(1_57_64_23_15_76_42_3);}_{(1_64_53_72_16_45_37_2);} & ^{4^3:}_{5^3:} & ^{(1_45_27_36_14_52_73_6);}_{(1_36_74_25_13_67_42_5);} \\
^{6^1:}_{7^1:} & ^{(1_72_36_54_17_23_65_4);}_{(1_76_32_45_17_63_24_5);} & ^{6^2:}_{7^2:} & ^{(1_67_52_43_16_75_24_3);}_{(1_27_46_53_12_74_65_3);} & ^{6^3:}_{7^3:} & ^{(1_26_47_35_12_64_73_5);}_{(1_62_57_34_16_25_73_4).}
\end{array}$$

\noindent Each 2-arc of $\G$ is suggested exactly once in these oriented cycles $i^j$. Each 2-path of $\G$ is suggested twice in them, once for each one of its two composing O-O 2-arcs.
The assumed orientation of each $\vec{C}_7^2=i_j$ corresponds with, and is induced by, the orientation of the corresponding 7-cycle $\vec{C}_7=\underline{i_j}$.\bigskip

\noindent Each 2-path $E$ of $\G$ separates two of its 24 7-cycles, say $\underline{i^j}$ and $\underline{k^\ell}$, with opposite orientations over $E$. Now, these $\underline{i^j}$ and $\underline{k^\ell}$ restrict to the two different 2-arcs provided by $E$, say 2-arcs $\vec{E}$ and $(\vec{E})^{-1}$.
Then, $\vec{E}$ and $(\vec{E})^{-1}$ represent corresponding arcs $\vec{e}$ and $(\vec{e})^{-1}$ in $i^j$ and $k^\ell$, respectively.\bigskip

\noindent Let us see that $\vec{e}$ and $(\vec{e})^{-1}$ can be zipped into an edge $e$ of $\G'$. In fact, $\G'$ can be assembled from the three charts shown on the upper right and bottom of Figure 1 by zipping the oriented 7-cycles $i^j$, interpreted all with counterclockwise orientation. Each of these three charts conforms a `rosette', where the oriented 7-cycles $i^j$ with $i\ne 0$ are represented as `petals' of the `central' oriented 7-cycles $0^1$, $0^2$ and $0^3$. (Similarly, the assembly of $\G'$ could have been done also around $i^1$, $i^2$ and $i^3$, taken as `central' oriented 7-cycles, for any $0\ne i\in J_7$).\bigskip

\noindent Moreover, each arc $\vec{e}$ in the external border of any selected one of the three charts, ($\vec{e}$ interpreted as an arc of an oriented cycle $\vec{C}_7^2$ in the selected chart), is accompanied, externally to the chart, by the symbol $i^j$ of another oriented 7-cycle $i^j$ that also contains $\vec{e}$ and forms a `petal' in just one of the other two (`rosette') charts. For example, the oriented cycle $3^1$ on the left of the chart centered at the oriented 7-cycle $0^1$ (on the upper-right of Figure 1) has its leftmost arc $\vec{e}$, corresponding to the symbol subsequence $6_52$ in $3^1=(1_34_76_52_13_47_65_2)$, also present in reverse in  the oriented cycle $1^3=(1_47_25_63_14_72_56_3)$, that is to say as $(\vec{e})^{-1}$, corresponding to the symbol subsequence $2_56$, at the upper-right in the chart centered at the oriented 7-cycle $0^3$ (on the lower-right of Figure 1). Thus, the symbols $1_3$ and $3_1$ accompany the representation of the arcs $\vec{e}$ and $(\vec{e})^{-1}$ on the outside of the external borders of their respective charts.
Not only the symbol $3^1$ indicates externally the arc $(\vec{e})^{-1}$ of the oriented 7-cycle $1^3$, but also indicates an arc $\vec{f}$ of the oriented 7-cycle $5^3$, the one corresponding to the symbol subsequence $6_74$ in $5^3=(1_36_74_25_13_67_42_5)$. The arc $(\vec{f})^{-1}$ is in the first mentioned oriented 7-cycle, $3^1$, just up from $\vec{e}$ and preceding it in the 28-cycle delimiting externally the chart centered at the 7-cycle $0^1$,
with corresponding symbol subsequence $4_76$ in $3^1=(1_34_76_52_13_47_65_2)$.\bigskip

\noindent The  presence of these arcs, $\vec{e}$, $(\vec{e})^{-1}$, $\vec{f}$ and $(\vec{f})^{-1}$, (and in all other similar cases) is expressed in the following formulation of the three oriented 28-cycles delimiting externally the charts with central oriented 7-cycles $0^1$, $0^2$ and $0^3$, (depicted respectively in the upper-right, lower-left and lower right of Figure 1), in the same color notation of the 24 oriented 7-cycles $i^j$ given above:
$$\begin{array}{c}
_{(1^1(3^36^3)\!(5^2)2^1(4^37^3)\!(6^2)3^1(5^31^3)\!(7^2)4^1(6^32^3)\!(1^2)5^1(7^33^3)\!(2^2)6^1(1^34^3)\!(3^2)7^1(2^35^3)\!(4^2))}
^{(4\,_6\,7\,_3\,2\,_1\,5\,_4\,6\,_1\,2\,_5\,4\,_3\,7\,_6\,1\,_3\,4\,_7\,6\,_5\,2\,_1\,3\,_5\,6\,_2\,1\,_7\,4\,_3\,5\,_7\,1\,_4\,3\,_2\,6\,_5\,7\,_2\,3\,_6\,5\,_4\,1\,_7\,2\,_4\,5\,_1\,7\,_6\,3\,_2)}\\
\\ _{(1^2(5^14^1)\!(2^3)3^2(7^16^1)\!(4^3)5^2(2^11^1)\!(6^3)7^2(4^13^1)\!(1^3)2^2(6^15^1)\!(3^3)4^2(1^17^1)\!(5^3)6^2(3^12^1)\!(7^3))}
^{(7\,_5\,2\,_6\,1\,_7\,5\,_3\,4\,_2\,6\,_3\,5\,_4\,2\,_7\,1\,_6\,3\,_7\,2\,_1\,6\,_4\,5\,_3\,7\,_4\,6\,_5\,3\,_1\,2\,_7\,4\,_1\,3\,_2\,7\,_5\,6\,_4\,1\,_5\,7\,_6\,4\,_2\,3\,_1\,5\,_2\,4\,_6\,1\,_6)}\\ \\
_{(1^3(2^27^2)\!(3^1)5^3(6^24^2)\!(7^1)2^3(3^21^2)\!(4^1)6^3(7^25^2)\!(1^1)3^3(4^22^2)\!(5^1)7^3(1^26^2)\!(2^1)4^3(5^23^2)\!(6^1))}
^{(5\,_6\,3\,_1\,4\,_7\,2\,_5\,6\,_7\,4\,_2\,5\,_1\,3\,_6\,7\,_1\,5\,_3\,6\,_2\,4\,_7\,1\,_2\,6\,_4\,7\,_3\,5\,_1\,2\,_3\,7\,_5\,1\,_4\,6\,_2\,3\,_4\,1\,_6\,2\,_5\,7\,_3\,4\,_5\,2\,_7\,3\,_6\,1\,_4)}
\end{array}$$

\noindent accompanying, below the part of each of these three oriented 28-cycle common with an oriented 7-cycle $i^j$, (like the initial $4_67_32_15_4\ldots$), with an integrated expression $i^j(\ldots,\ldots)(\ldots)$ containing, between the first pair of parentheses, $(\ldots,\ldots)$, the symbols of the oriented 7-cycles containing $\vec{f}$ and $(\vec{e})^{-1}$ in the other two charts in each case, where $(\vec{f})^{-1}$ and $\vec{e}$ are the corresponding arcs in $i^j$, and containing, between the second pair of parenthesis, $(\ldots)$, the symbol following them externally to the chart involved, in counterclockwise fashion, (like the immediately lower accompanying $1^1(3^36^3)(5^2)\ldots$).\bigskip

\noindent This codifies the assembly of the three charts into the claimed graph $\G'$. Moreover,  the $24$ oriented 7-cycles $i^j$ can be filled each with a corresponding 2-cell, so that because of the cancelations of the two opposite arcs on each edge of $\G'$ (for having opposite orientations makes them mutually cancelable), $\G'$ becomes embedded into a closed orientable surface $T_3$. As for the genus of $T_3$, observe that $$|V(\G')|=2\times 28=56\,\,\mbox{ and }\,\,|E(\G')|=2|E(\G)|=2\times 42=84,$$
so that by the Euler characteristic formula for $T_3$ here,
$$|V(\G')|-|E(\G')|+|F(\G')|=56-84+24=-4=2-2.g(T_3),$$ and thus $g=3$, so $T_3$ is a 3-torus.
This yields the Klein map of Coxeter notation $\{7,3\}_8$. (See \cite{Sch,K,SL} and note that the Petrie polygons of this map are 8-cycles).\bigskip

\begin{thm} The Klein graph $\G'$ is both a $\{C_7\}_{P_2}$-UH graph and a $\{\vec{C}_7\}_{\vec{P}_2}$-UH digraph, composed by $24$  {\rm(}oriented{\rm)} $7$-cycles that yield the Klein map $\{7,3\}_8$ in $T_3$.\qfd\end{thm}

\noindent For the Klein map $\{7,3\}_8$, the 3-torus appeared originally dressed as the Klein quartic $x^3y+y^3z+z^3x=0$, a Riemann surface and the most symmetrical curve of genus 3 over the complex numbers. The automorphism group for this Klein map is $PSL(2,7)=GL(3,2)$, (\cite{BL}), the same automorphism group of $\mathcal F$, whose index is 2 in the common automorphism groups of $\G$, $\G'$ and $\G''$.\bigskip

\begin{cor} The Klein quartic graph ${\mathcal K}$, whose vertices are the $7$-cycles $i^j$ of $\G'$, with adjacency between two vertices if their representative $7$-cycles have a pair of {\rm O-O} arcs, is regular of degree $7$, chromatic number $8$ and has a natural triangular $T_3$-embedding yielding the dual Klein map $\{3,7\}_8$.\end{cor}\bigskip

\proof Each vertex $i^j$ of ${\mathcal K}$ is assigned color $i\in\{0\}\cup J_7$. Also, we have a partition of $T_3$ into 24 connected regions, each region having exactly seven neighboring regions, with eight colors needed for a proper map coloring. \qfd

\section{Final remarks}

Following the remarks made after Theorems 1 and 2, it can be said that the zipping method of Section 3 can be adapted to other graphical situations; to begin with, the Pappus graph, the Desargues graph and the Biggs-Smith graph, the last one yielding
the Men\-ger graph of a self-dual $(102_4)$-con\-fi\-gu\-ra\-tion,
what may be called a $\{K_4,L(Q_3)\}_{K_3}$-UH graph, in a similar way in which the graph of \cite{D1} is a $\{K_4,K_{2,2,2}\}_{K_2}$-UH graph, where $L(Q_3)$ is the line graph of the 3-cube  graph $Q_3$. More specifically, the Biggs-Smith graph yields, by means of an adequate zipping procedure, a connected 12-regular graph which is the union of 102 copies of $L(Q_3)$ without common squares as well as the edge-disjoint union of 102 copies of $K_4$, with each triangle (edge) as the intersection of exactly
two (four) copies of $L(Q_3)$.
Also, generalizing on zipping results over the Desargues graph, it can be concluded that the line graph $L(K_n)$, with $n\geq 4$, is a tightly fastened $\{K_{n-1},K_3\}_{K_2}$-UH graph with $n$ copies of $K_{n-1}$ and ${n\choose 3}$ copies of $K_3$.

A final remark is that the role played by the Heawood graph $\G''$ in the construction of the so-called Ljubljana semi-symmetric graph \cite{BDT,Cetal}, which is an 8-cover of $\G''$, makes us wonder whether there are any more relations between this 8-cover and both $\G$ and $\G'$, derived all ultimately from $\G''$.


\begin{thebibliography}{99}
\bibitem{Bondy} J. A. Bondy, {\it Variations of the hamiltonian theme}, Canad. Math. Bull., {\bf 15} (1972), 57--62.
\bibitem{F} I. Z. Bouwer et al., The Foster Census, R. M. Foster's Census of Connected Symmetric Trivalent Graphs, Charles Babbage Res. Ctr., Canada 1988.
\bibitem{BDT} A. E. Brouwer, A. J. Dejter and C. Thomassen, {\it Highly symmetric subgraphs of hypercubes}, J. Algebraic Combin.,
     {\bf 2}(1993) 25--29.
\bibitem{BL} E. Brown and N. Loehr, {\it Why is $PSL(2, 7)=GL(3, 2)?$}, Amer. Math. Mo., {\bf 116-8}, Oct. 2009,
727--732.
\bibitem{Cam} P. J. Cameron, {\it $6$-transitive graphs}, J. Combin. Theory Ser. B {\bf 28}(1980), 168-179.
\bibitem{Cher} G. L. Cherlin, The Classification of Countable Homogeneous Directed Graphs and Countable Homogeneous $n$-tournaments, Memoirs Amer. Math. Soc., vol. 131, number 612, Providence RI, January 1988.

\bibitem{Cetal} M. Conder, A. Malnic, D. Marušic, T. Pisanski, and P Potocnik, {\it The edge-transitive but not vertex-transitive cubic graph on 112 vertices}, Jour. Graph Theory,
{\bf 50} (2005), 25-42.

\bibitem{Cox} H. S. M. Coxeter, {\it Self-dual configurations and regular graphs},
Bull. Amer. Math. Soc., {\bf 56}(1950), 413--455.

\bibitem{D1} I. J. Dejter, {\it On a $\{K_4,K_{2,2,2}\}$-ul\-tra\-ho\-mo\-ge\-neous graph},
Australasian Journal of Combinatorics, {\bf 44}(2009), 63--76.

\bibitem{Dx} I. J. Dejter, {\it On a $\vec{C}_4$-ultrahomogeneous oriented graph}, Discrete Mathematics, {\bf 310}(2010), 1389--1391.
\bibitem{FH} G. Fan and R. Haggkvist, {\it The Square of a Hamiltonian Cycle},
SIAM Jour. Discrete Math., {\bf 7} (1994), 203--212.
\bibitem{Fra} R. Fra\"\i  ss\'e, {\it Sur l'extension aux relations de quelques propriet\'es des ordres}, Ann.
Sci. \'Ecole Norm. Sup. 71 (1954), 363--388.
\bibitem{Gard} A. Gardiner, {\it Homogeneous graphs}, J. Combinatorial Theory (B), {\bf 20}
(1976), 94--102.
\bibitem{GR} C. Godsil and G. Royle, Algebraic Graph Theory, Springer--Verlag, 2001.
\bibitem{GK} Ja. Ju. Gol'fand and M. H. Klin, {\it On $k$-homogeneous graphs}, Algorithmic studies in combinatorics (Russian), {\bf 186}(1978), 76--85.
\bibitem{I} D. C. Isaksen, C. Jankowski and S. Proctor, {\it On $K_*$-ul\-tra\-ho\-mo\-ge\-neous graphs}, Ars Combinatoria, Volume LXXXII, (2007), 83--96.
\bibitem{Lach} A. H. Lachlan and R. Woodrow, {\it Countable ultrahomogeneous undirected
graphs}, Trans. Amer. Math. Soc. 262 (1980), 51–-94.
\bibitem{SL} S. Levy (ed.), The Eightfold Way: The Beauty of the Klein Quartic,  Cambridge University Press, New York, 1999.
\bibitem{K}
F. Klein, {\it \"Uber die Transformationen siebenter Ordnung der elliptischen Funktionen},  Math. Ann., {\it 14} (1879), 428--471, 1879.

\bibitem{Ronse} C. Ronse, {\it On homogeneous graphs}, J. London Math. Soc., {\bf s2-17} (1978), 375--379.
\bibitem{Sheh} J. Sheehan, {\it Smoothly embeddable subgraphs}, J. London Math. Soc.,
{\bf s2-9} (1974), 212--218.
\bibitem{Sch} E. Schulte and J. M. Wills, {\it A Polyhedral Realization of Felix Klein's Map $\{3, 7\}_8$ on a Riemann Surface of Genus 3},
J. London Math. Soc., {\bf s2-32} (1985), 539--547.

\end{thebibliography}
\end{document}